\documentclass[reqno,11pt]{amsart} 

\usepackage{amsmath, amsthm, a4, latexsym, amssymb, pstricks, natbib, setspace}
\usepackage[dvips]{graphicx}
\usepackage{natbib}

\def\@rmrk#1#2{\refstepcounter
    {#1}\@ifnextchar[{\@yrmrk{#1}{#2}}{\@xrmrk{#1}{#2}}}

\makeatletter\@addtoreset{equation}{section}\makeatother

 \newfont{\bfit}{cmbxti10 scaled 1200}

 \newcommand{\eps}{\varepsilon}

 \newcommand{\R}{\mathbb{R}}
 \newcommand{\N}{\mathbb{N}}

 \newcommand{\prob}{\mathbb{P}}

 \newcommand{\me}{\mathbb{E}}
 
 \newcommand{\1}{{\sf 1}}

 \newcommand{\sfrac}[2]{\mbox{$\frac{#1}{#2}$}}
 
 \newcommand{\ssup}[1] {{\scriptscriptstyle{({#1}})}}

\renewcommand{\subsection}{\secdef \subsct\sbsect}
\newcommand{\subsct}[2][default]{\refstepcounter{subsection}
\vspace{0.15cm}
{\flushleft\bf \arabic{section}.\arabic{subsection}~\bf #1  }
\nopagebreak\nopagebreak}
\newcommand{\sbsect}[1]{\vspace{0.1cm}\noindent
{\bf #1}\vspace{0.1cm}}

{\nopagebreak {\hfill\rule{2mm}{2mm}}\\ }

\newtheorem{theorem}{Theorem}[section]
\newtheorem{thm}{Theorem}
\newtheorem{lemma}[theorem]{Lemma}

\newcommand{\p}{\mathbb{P}}
\newcommand{\ra}{\rightarrow}
\newcommand{\E}{\mathbb{E}}
\newcommand{\ta}{\tilde{a}}
\newcommand{\fnm}{\tfrac{\nu}{\mu}}
\newcommand{\bfnm}{(\tfrac{\nu}{\mu})}

\newtheoremstyle{thm}{1.5ex}{1.5ex}{\itshape\rmfamily}{}
{\bfseries\rmfamily}{}{2ex}{}

\newtheoremstyle{rem}{1.3ex}{1.3ex}{\rmfamily}{}
{\itshape\rmfamily}{}{1.5ex}{}
\theoremstyle{rem}

\refstepcounter{subsubsection}

\begin{document}

\title[Small value probabilities]
{Small value probabilities via the branching tree heuristic}

\centerline{\LARGE \bf Small value probabilities via the branching tree heuristic}

\vspace{0.5cm}

\thispagestyle{empty}
\vspace{0.2cm}
\textsc{Peter M\"orters\footnote{Communicating author.}} and \textsc{Marcel Ortgiese}\\
Department of Mathematical Sciences,
University of Bath, Bath BA2 7AY, England\\
E--mail: \texttt{maspm@bath.ac.uk} and \texttt{ma2mo@bath.ac.uk}

\vspace{0.3cm}


\begin{quote}{\small {\bf Abstract: } In the first part of this paper we give easy and intuitive
proofs for the small value probabilities of the martingale limit of a supercritical Galton-Watson process
in both the Schr\"oder and the B\"ottcher case. These results are well-known, but the most cited
proofs rely on generating function arguments which are hard to transfer to other settings. 
In the second part we show that the strategy underlying our proofs can be used in the quite 
different context of self-intersections of stochastic processes. Solving a problem posed by Wenbo Li, 
we find the small value probabilities for intersection local times of several Brownian motions, as well 
as for self-intersection local times of a single Brownian motion.}
\end{quote}
\vspace{0.5cm}

{\footnotesize
{\bf MSc classification (2000):} Primary 60F10 Secondary  60J80, 60J65, 60J55.\\[-5mm] 

}



\section{Introduction}

The \emph{small value problem} is to find, for a nonnegative random variable~$X$,
the speed of decay of the left tail $\prob\{X <\eps\}$ as $\eps\downarrow 0$. 
Important examples are the \emph{small ball problem} where $X$ is the norm of
a random variable with values in a Banach space, the \emph{lower level problem}
where $X$ is the maximum of a continuous random process $(X(t) \colon t\in[0,1])$,
or \emph{boundary crossing problems} where $X$ is the first exit
time of a stochastic process from a general space-time domain.\smallskip

Small value problems arise in a great variety of contexts in probability
and analysis. Examples include approximation and quantisation 
problems \citep{LL99, DF+03, GLP03}, Brownian pursuit problems \citep{LS01b},
polymer measures \citep{HHK97}, and convex geometry \citep{KV07}.
A systematic theory of small value problems, however, is only available 
when $X$ is the norm of a Gaussian random variable. 
For other cases some isolated
techniques are known, but a bigger picture has not yet emerged. A survey of Gaussian methods
in this field is \cite{LS01a} and an updated bibliography on small value problems is kept 
at \cite{Li06}.\smallskip

In this paper we contribute to the theory of small value problems by presenting
systematically an approach which we found successful in a variety of cases. We illustrate
our technique by three main examples. The \emph{first} example is the most natural one for our approach:
the martingale limit of a supercritical Galton-Watson process. In this case the small value problem
has been solved --- by~\cite{Du71b, Du71a} in the Schr\"oder case and, up to a Tauberian theorem
of~\cite{Bi88}, also in the B\"ottcher case. These proofs use an integral transformation
approach together with some nontrivial complex analysis, a powerful method, but inflexible
and not very intuitive. Our method, by contrast, is very simple and based on an easy intuition. {F}rom this
example we derive the term \emph{`branching tree heuristic'} for the general approach.
\smallskip

\pagebreak
The \emph{second} example is our main result and treated here for the first time: 
We solve a problem posed by Wenbo Li at the
Miniworkshop `Small deviation probabilities and related topics' at Oberwolfach in October~2003.
The problem is to identify the small value probability of the random variable
$$X=\int_{-\infty}^\infty \prod_{i=1}^m L_i^{q_i}(x,1)\, dx\, ,$$
where $L_1(x,t),\ldots, L_m(x,t)$ are the local times of $m\ge 2$ independent Brownian motions. We explain very
carefully how a heuristic embedding of a tree in the Brownian motion framework leads to a proof
based on the same principles as in the Schr\"oder case of the first example.\smallskip

Also our \emph{third} example appears to be new, though it is really quite elementary. 
We look at the $L^q$-norm of the local time of
a single Brownian motion stopped when it exits a bounded interval for the first time, which, for 
$q$ an integer, may be interpreted as the $q$-fold self-intersection local time of the motion.
We again find a relation to a Galton-Watson tree, this time of B\"ottcher type, and exploit
this relation to find a strikingly simple proof of the small value probability.\smallskip
 
We believe that our method can be used in a number of further cases, when the optimal strategy
for a random variable to obtain small values is inhomogeneous. We conclude the paper
with an outlook to future research.\medskip

\section{Small value probabilities for the martingale limit of a Galton-Watson tree}

Consider a Galton-Watson branching process $(Z_n \colon n\geq 0)$ with offspring distribution
$(p_k \colon k\geq 0)$ starting with a single founding ancestor, called $\rho$, in generation $0$. We suppose
that the offspring variable $N$ is nondegenerate and satisfies $\mu := \E N > 1$ and $\E [ N \log N ] < \infty$.
By the famous Kesten-Stigum theorem these conditions ensure that the martingale limit
$$W:= \lim_{n\to\infty} \frac{Z_n}{\mu^n}$$
exists and is nontrivial almost surely on survival. Except in the case when $N$ is
geometric, the distribution of $W$ is not known explicitly and one relies on asymptotic results to
describe its behaviour.

For the formulation of our results, we further assume $p_0 = 0$, which is no
loss of generality: Removing all finite subtrees from a Galton-Watson tree does not
change its martingale limit, but the resulting tree is still a Galton-Watson tree (with a modified
offspring variable), see \citet[Chapter~1, Section~12]{AN72}. 

As usual we distinguish between the \emph{Schr\"oder} case and the \emph{B\"ottcher} case,
depending on whether $p_1 > 0$ or $p_1 = 0$. These two cases yield very different lower
tail behaviour for $W$. In the following $a(\eps) \asymp b(\eps)$ means that there exist
constants $0<c<C<\infty$ such that $$c a(\eps) \le b(\eps) \le C a(\eps),
\qquad\mbox{ for all $0<\eps<1$.}\\[3mm]$$

\begin{thm}[Dubuc 1971]\ \\ \vspace{-0.5cm} \label{main}
\begin{itemize}
\item[(a)] In the Schr\"oder case define $\tau:=-\log p_1/\log \mu>0$. Then
\[ \p\{ W < \eps \} \asymp \eps^\tau \, .\] \\[-5mm]
\item[(b)] In the B\"ottcher case define $\nu := \min\{ i \geq 0 \ : \ p_i \neq 0\}\ge 2$
and $\beta := \frac{\log \nu}{\log \mu} < 1$. Then
$$-\log \p\{ W < \eps \} \asymp \eps^{\frac{-\beta}{1-\beta}} \, .$$\\[-10mm]
\end{itemize}
\end{thm}

In this paper we offer simple proofs of both parts of Theorem~\ref{main}, and show how the idea
behind these proofs can be adapted to obtain small value probabilities for situations,
which might look quite different at a first glance.

The main idea of the proofs is to understand the optimal strategy by which the tree keeps the generation
size small. It turns out that the best strategy consists of producing as little offspring as possible at the
beginning and then, once the necessary reduction in size is achieved, letting the tree grow normally.
If the tree produces a larger number of children at the beginning, it will be more expensive to
control the growth later on, since every additional child is likely to produce more than one child as well.
This effect is illustrated in Figure~1.

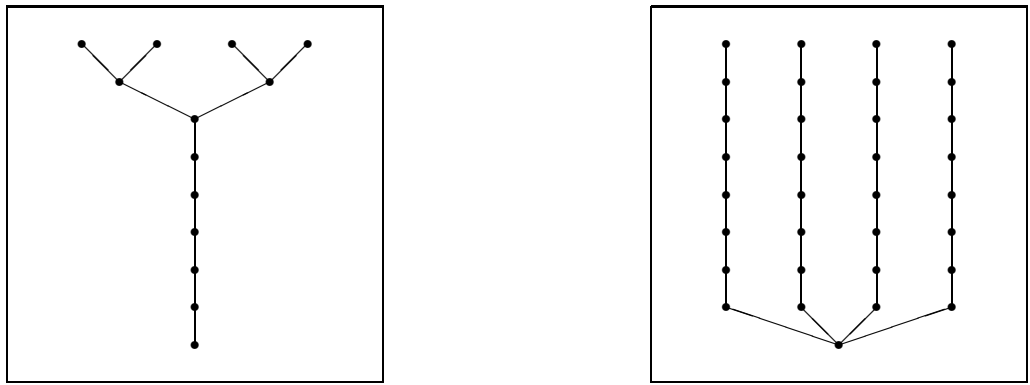
\begin{figure}[htbp]
	\hfill
	\begin{minipage}[t]{.45\textwidth} 
		\begin{center}
		\setlength{\unitlength}{0.5cm}
			\begin{picture}(10,10)
			\put(0,0){\line(0,10){10}}
			\put(0,0){\line(10,0){10}}
			\put(0,10){\line(10,0){10}}
			\put(10,10){\line(0,-10){10}}
			\multiput(5,1)(0,1){7}{\circle*{0.2}}
			\put(5,1){\line(0,1){6}}

			\put(3,8){\circle*{0.2}}
			\put(7,8){\circle*{0.2}}
			\put(5,7){\line(-2,1){2}}
			\put(5,7){\line(2,1){2}}

			\multiput(2,9)(2,0){4}{\circle*{0.2}}
			\put(3,8){\line(-1,1){1}}
			\put(3,8){\line(1,1){1}}
			\put(7,8){\line(-1,1){1}}
			\put(7,8){\line(1,1){1}}
			\end{picture}
		\end{center}
	\end{minipage} 
	\hfill
	\begin{minipage}[t]{.45\textwidth}
		\begin{center}
		\setlength{\unitlength}{0.5cm}
			\begin{picture}(10,10)
			\put(0,0){\line(0,10){10}}
			\put(0,0){\line(10,0){10}}
			\put(0,10){\line(10,0){10}}
			\put(10,10){\line(0,-10){10}}
			\put(5,1){\circle*{0.2}}
			\multiput(2,2)(0,1){8}{\circle*{0.2}}
			\multiput(4,2)(0,1){8}{\circle*{0.2}}
			\multiput(6,2)(0,1){8}{\circle*{0.2}}
			\multiput(8,2)(0,1){8}{\circle*{0.2}}
			\put(5,1){\line(-3,1){3}}
			\put(5,1){\line(-1,1){1}}
			\put(5,1){\line(1,1){1}}
			\put(5,1){\line(3,1){3}}
			\put(2,2){\line(0,1){7}}
			\put(4,2){\line(0,1){7}}
			\put(6,2){\line(0,1){7}}
			\put(8,2){\line(0,1){7}}
			\end{picture}
		\end{center}
	\end{minipage}
	\caption{The picture on the left illustrates the optimal strategy to keep the final generation size small. 
	By comparison in the picture on the right the offspring of more individuals have to be kept under control to   
	produce the same effect.}
\end{figure}

By $(Z_n(v) \colon n\geq 0)$ we denote the generation sizes of the subtree consisting of all
the descendants of the individual $v$. Note  that for each fixed~$v$ the process
$(Z_n(v) \colon n\geq 0)$ is again a Galton-Watson process and hence we can define 
the martingale limit 
\[ W(v) := \lim_{n\ra \infty} \frac{Z_n(v)}{\mu^n}. \]
Let $v_k(1), \ldots, v_k(Z_k)$ be the individuals in the $k^{\rm th}$ generation. By decomposing the 
individuals in the $n^{\rm th}$ generation according to their ancestors in the $k^{\rm th}$
generation we get, for all $n\ge k$,
$$Z_n= \sum_{i=1}^{Z_k} Z_{n-k}(v_k(i))\, .$$
Hence we obtain
\begin{equation} 
W = \lim_{n\ra \infty} \frac{Z_n}{\mu^n} =  \lim_{n\ra\infty} \mu^{-k}\, 
\sum_{i=1}^{Z_k} \frac{Z_{n-k}(v_k(i))}{\mu^{n-k}}
= \mu^{-k} \sum_{i=1}^{Z_k} W(v_k(i)) ,\label{eqn:branching_structure} \end{equation}
where all the random variables $W(v_k(i))$ are iid with the same distribution as $W$.\\[1mm]

This section is organised as follows: We first investigate the Schr\"oder case. We start by showing
that the suggested strategy is successful, which proves the lower bound. We then give a rough
argument which produces the precise logarithmic asymptotics. This argument is then
refined, exploiting the self-similarity of the tree, to complete the proof of Theorem~\ref{main}\,(a).
The arguments leading to the result in the B\"ottcher case, Theorem~\ref{main}\,(b), 
are easier and given in the final two subsections.

\pagebreak[3]

\subsection{The Schr\"oder case: The lower bound}\label{schroeder}

For the \emph{lower bound} suppose $0<\eps<1$ and pick $n$
such that $\mu^{-n}\le \eps< \mu^{-n+1}$.
Using (\ref{eqn:branching_structure}) we obtain
$$\begin{aligned} \p\{ W < \eps \} & \geq \p\{ W < \mu^{-n} \, | \, Z_n = 1\}\, \p\{Z_n = 1\}
\\ & = \p\big\{\mu^{-n} \, W(v_n(1)) < \mu^{-n} \big\}\, p_1^n = c \,  p_1^n \geq (cp_1) \eps^\tau\, , 
\end{aligned}$$
where $c := \p\{W < 1\}>0$.

\subsection{The Schr\"oder case: The logarithmic upper bound}

As the \emph{first step} in the proof of the upper bound we show that
\begin{equation}\label{roughbound}
\limsup_{\eps\downarrow 0} \frac{\log \prob\{ W< \eps \}}{-\log \eps} \le  - \tau\, .
\end{equation}

{\bf Remark:} In the second step of the argument we only use that $\prob\{ W<\eps\}$ decreases
like \emph{some} positive power of $\eps$. Other instances of our method, however, make use of lower
bounds on this power, so it is instructive to show the `best possible' argument here.\qed\\

Fix a large $m$ for the moment, and let $n\ge m$.
By decomposing the set of individuals in the $n^{\rm th}$~generation of
the branching process according to their last common ancestor with the `spine'
$\rho=v_0(1), v_1(1), v_2(1), \ldots, v_m(1)$
consisting of the leftmost individual in each of the first $m+1$ generations, we obtain a decomposition
$$Z_n = \sum_{k=1}^m \sum_{j=2}^{Z_1(v_{k-1}(1))}   Z_{n-k}(v_k(j)) + Z_{n-m}(v_m(1))\, .$$
Discarding the contributions for $j\ge 3$, if they exist, and also the last summand, dividing by $\mu^n$ and
letting $n\uparrow\infty$, gives
\begin{equation}\label{spinaldecomp}
W \ge \sum_{k=1}^m \mu^{-k} W_k\, ,
\end{equation}
where $W_k=0$ if $v_{k-1}(1)$ has only one offspring, and $W_k=W(v_{k}(2))$ otherwise.
Note that $W_1, \ldots, W_k$ are independent, identically distributed with distribution given by
$\prob\big\{ W_k=0 \} =p_1$ and $$\prob\big\{ W_k<x \, | \, W_k \not= 0 \} = \prob\{ W<x \}
\qquad \mbox{ for all } x>0 \, .$$

\begin{figure}[htbp]
\centering
\includegraphics[width=8cm, height=5cm]{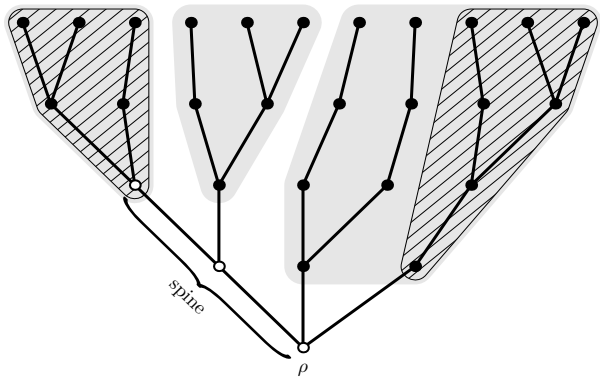} 
\caption{Decomposition of the tree according to the ancestry from a spine with length $m=2$. 
The shaded parts of the tree are discarded in our calculation.}
\end{figure}

Now suppose $\delta>0$ is given. As $W>0$ almost surely, there exists $\theta>0$ such that
$\prob\{W<\theta\} \le \delta p_1$. We fix the integer $\ell$ such that $\mu^{\ell}\le \theta<\mu^{\ell+1}$.
Let $\eps>0$ be arbitrary and define $n$ by $\mu^{-n-1}< \eps \le \mu^{-n}$.
Then, using \eqref{spinaldecomp} for $m=n+\ell$,
$$\begin{aligned}
\prob\big\{ W<\eps \big\}  & \le \prob\big\{ W<\mu^{-n} \big\}
\le \prob\Big\{ \sum_{k=1}^m \mu^{-k} W_k <\mu^{-n} \Big\}
\le \prod_{k=1}^m \prob\big\{ W_k <\mu^{-n+m} \big\} \\
& \le \big( p_1 + \prob\big\{ W <\theta \big\} \big)^m
\le \Big( p_1^\ell (1+\delta)^\ell\Big) \, p_1^{n} \, (1+\delta)^n
\le C\, \eps^\tau \, e^{\delta n}\, ,
\end{aligned}$$
for $C:=  p_1^\ell \,(1+\delta)^\ell \, \mu^\tau$, from which \eqref{roughbound} follows, as $\delta>0$ was arbitrary.

\subsection{The Schr\"oder case: Up-to-constants asymptotics}

We are now in a position to refine the upper bound and prove Theorem~\ref{main}\,(a).
Define a sequence $(a(n) \colon n\ge 0)$ by setting 
$$a(n) := \p\{W < \mu^{-n}\} \,p_1^{-n}.\\[1mm]$$
For arbitrary $0<\eps<1$ we pick the integer $n\ge 0$ such that $\mu^{-n-1}\le \eps < \mu^{-n}$.
Then\\[-2mm] $$\prob\big\{ W < \eps \big\} \le \p\{ W < \mu^{-n} \}
= a(n)\,p_1^n \le a(n) \, (1/p_1)\, \eps^\tau\, , \\[1mm]$$
hence, to complete the proof, it suffices to show that $(a(n) \colon n\ge 0)$ is bounded.

Denote by $N_n$ the number of offspring of the left-most individual in \mbox{generation $n$},
and let $$T := \min\{ n \geq 0 \, : \, N_n \neq 1\}\, .$$
Obviously, \mbox{$\p\{ T = j\} = p_1^{j} (1-p_1)$}. Let $j < n$ be nonnegative integers.
Applying (\ref{eqn:branching_structure}) we get
\begin{equation}\begin{aligned} \p \big\{ W < \mu^{-n},  \,  T = j \big\}
&  \le  \p\big\{ \mu^{-(j+1)} \big( W(v_{j+1}(1)) + W(v_{j+1}(2)) \big) < \mu^{-n} , T = j \big\} \\[2mm]
& \le p_1^{j+1} \, \prob\big\{ W< \mu^{-(n-j-1))} \big\}\, \beta(n-j-1) \, ,\label{estimate_for_W_fixed_T}
\end{aligned}\end{equation}
where $\beta(i):= p_1^{-1} \, \p\{W < \mu^{-i}\}$.
By the a-priori estimate \eqref{roughbound} we have $\sum \beta(i) < \infty$.

Using (\ref{estimate_for_W_fixed_T}) we get, for any positive integer~$n$,
\begin{equation}\begin{aligned}
\p\{ W < \mu^{-n} \} &  
\leq \ \sum_{j=0}^{n-1} \p\{ W < \mu^{-n} ,\, T = j\} + \p\{T \ge n\} \\
& \leq \sum_{j=0}^{n-1} \,p_1^{j+1} \, \p\{W < \mu^{-(n-j-1)}\} \, \beta(n-j-1) + p_1^n \, .
\label{inequality_for_W}\end{aligned}\end{equation}

We deduce from~(\ref{inequality_for_W}) that
$a(n)  \leq  \sum_{j=0}^{n-1} \, a(n-j-1) \, \beta(n-j-1) + 1 \, .$
Define $\tilde{a}(-1) := 1$, $\beta(-1):=1$, and inductively, for nonnegative $n$,
$$ \ta(n)  := \sum_{j=0}^{n-1} \, \tilde a(n-j-1) \, \beta(n-j-1) + 1
= \sum_{j=-1}^{n-1} \, \ta(j) \, \beta(j) \, .$$
Then, since $a(n) \leq \ta(n)$ for all $n\ge 0$, it suffices to show that $(\ta(n) \colon n \ge 0)$ 
is bounded. {F}rom the definition it follows easily that $\ta(n)  = \ta(n-1)\, ( 1 + \beta(n-1) \,)$, hence
$\ta(n) = \prod_{i=0}^{n-1} \big( 1+\beta(i) \big) \, , $
which converges as $\sum_{i=0}^{\infty} \beta(i)$ converges. Hence $(\ta(n) \colon n\ge 0)$ is bounded
and the proof complete.


\subsection{The B\"ottcher case: The lower bound}\label{boettcherlow}

We now consider the case when $p_1 = 0$. Recall that $\nu := \min\{ j \geq 0 \, : \, p_j \neq 0\} \geq 2$
and $\nu<\mu$. For every $n$, there are at least $\nu^n$ individuals in generation $n$, hence
$$ \p\{Z_n = \nu^n \} = \p \{Z_n = \nu^n \,| \, Z_{n-1} = \nu^{n-1} \} \, \p \{ Z_{n-1} = \nu^{n-1} \} = 
p_\nu^{\nu^{n-1}}  \p \{ Z_{n-1} = \nu^{n-1} \} \, . $$
Also $\p \{ Z_1 = \nu\} = p_\nu$, and therefore
\begin{equation} \p\{Z_n = \nu^n \}  = p_\nu^{1 + \nu + \cdots + \nu^{n-1}} = p_\nu^{\frac{\nu^n -1}{\nu -1}}. 
\\[1mm] \label{eqn:Z_n_minimal}
\end{equation}
Given~$\eps>0$ we look at the \emph{lower bound} of the probability~$\prob\{W<\eps\}$. Pick the integer~$n$ 
such that $(\sfrac\nu\mu)^n \le \eps < (\sfrac\nu\mu)^{n-1}$. Invoking (\ref{eqn:branching_structure}) and (\ref{eqn:Z_n_minimal}) we get
\begin{eqnarray*}
\p\left\{ W < \eps \right\} & \geq & \p \big\{ W < \bfnm^n \,\big| \, Z_{n+1} = \nu^{n+1}\big\} \,
\p\left\{ Z_{n+1} = \nu^{n+1} \right\} \\
    & = & \p\big\{ W(v_{n+1}(1)) + \cdots + W(v_{n+1}(\nu^{n+1})) <
\left(\tfrac{\mu}{\nu}\right) \nu^{n+1} \big\}
\, p_\nu^{\frac{\nu^{n+1}-1}{\nu - 1}}\\
 & \geq & \p\Big\{ \,\Big| \! \sum_{j=1}^{\phantom{1}\nu^{n+1}}\!W(v_{n+1}(j)) - \nu^{n+1} \Big| < \delta \,
 \nu^{n+1}  \Big\} \, p_\nu^{\frac{\nu^{n+1}-1}{\nu - 1}} \ ,
\end{eqnarray*}
where $\delta := \sfrac{\mu}{\nu}-1 > 0$. By the weak law of large numbers we may choose $N \in \N$ such that
$$ \p\Big\{ \, \Big| \! \sum_{j=1}^{\phantom{1}\nu^{m+1}} W(v_{m+1}(j)) - \nu^{m+1} \Big|
< \delta \, \nu^{m+1} \Big\} \ge p_\nu^{1/(\nu -1)}\qquad \mbox { for all $m\ge N$.}$$
Then, for all $n \geq N$, we have
$$ - \log \p \left\{ W < \eps \right\} \leq (- \log p_\nu) \, {\frac{\nu^{n+1}}{\nu -1}}  \le C\eps^{\frac{-\beta}{1-\beta}}\, , $$
where $C := (- \log p_\nu) \frac{\nu^2}{\nu -1}$,
using that $(\sfrac\nu\mu)^{\frac{-\beta}{1-\beta}}=\nu$, by definition of $\beta$.

\subsection{The B\"ottcher case: The upper bound}\label{boettcherup}

Given $\eps>0$ we continue with an \emph{upper bound} for the probability~$\prob\{W<\eps\}$. 
Pick the integer $n$ such that $(\sfrac\nu\mu)^{n+1} \le \eps < (\sfrac\nu\mu)^{n}$. 
Using once again (\ref{eqn:branching_structure}) we get
\begin{equation}\label{eqn:link}
\p\left\{W < \eps \right\}
\le \p \Big\{ {\mu^{1-n}} \sum_{j=1}^{\nu^{n-1}} W(v_{n-1}(j)) <  \bfnm^n \Big\}
=  \p\big\{ S(\nu^{n-1}) > 0 \big\} \, ,
\end{equation}
where $X_j := \frac{\nu}{\mu} - W(v_{n-1}(j))$ and $S(k) := \sum_{j=1}^k X_j$. 

We now estimate the right hand side by a simple large deviation bound, which only uses that
$X_j$ is bounded from above and has negative mean. By the exponential Chebyshev inequality,
\begin{equation} \prob\big\{ S(k) \ge 0 \big\}\le
\p\{e^{\tau S(k)} \geq 1 \} \leq \E e^{\tau S(k)} = \left( \E
e^{\tau X_1}\right)^k \label{eqn:chebyshev}. \\[1mm]\end{equation} We claim there
exists $\tau > 0$ such that $\E e^{\tau X_1}< 1$. Indeed, denoting
$\varphi(\tau) := \E e^{\tau X_1}$ and using Lebesgue's dominated
convergence theorem, we have
$$\lim_{\tau \downarrow 0} \frac{\varphi(\tau) - \varphi(0)}{\tau}   =  \lim_{\tau \downarrow 0}
\E \left[\frac{ e^{\tau X_1} - 1}{\tau}\right] =  \me \lim_{\tau \downarrow 0} 
\left( \frac{e^{\tau X_1} -1}{\tau} \right)  = \E X_1 = \fnm - 1 < 0\, .$$
Since $\varphi(0) = 1$, we can thus choose $\tau > 0$ such that $\varphi(\tau)< 1$.
Combining this with (\ref{eqn:link}) and (\ref{eqn:chebyshev}), we get
$- \log \p\{ W < \eps \} \geq (- \log \varphi(\tau)) \, \nu^{n-1} \ge c \, \eps^{\frac{-\beta}{1-\beta}} \, , $
where $c := - \nu^{-2}\log \varphi(\tau) > 0$.

\section{Small value probabilities for mutual intersection local times}\label{mutual}

In this section we identify the small value probability of the random variables
$$X(t_1,\ldots, t_m):=\int_{-\infty}^\infty \prod_{i=1}^m L_i^{q_i}(x,t_i)\, dx\, ,$$
where $(L_1(x,t) \colon x\in\R, t\ge 0), \ldots, (L_m(x,t)\colon x\in\R, t\ge 0)$ 
are the local time fields of $m$ independent Brownian motions started at the origin.
For $q_1=\cdots=q_m=1$ the random variable $X(t_1,\ldots, t_m)$ measures
the amount of intersection between the motions up to times  
$t_1,\ldots, t_m$ and it is therefore called \emph{(mutual) intersection local time}.
\smallskip

Our solution to the small value problem for intersection local times
is based on an analogy between the martingale limit $W$ of a
Galton-Watson tree in the Schr\"oder case and the random variables $X(\sigma^{\ssup 1},\ldots,
\sigma^{\ssup m})$ where $\sigma^{\ssup 1},\ldots, \sigma^{\ssup m}$ are the first exit times of 
the Brownian motions from the interval~$(-1,1)$. This analogy allows us to carry over the crucial 
steps in the proof of Theorem~\ref{main}\,(a) to the new situation, and hence to prove the following theorem.

\begin{thm}\label{main2}
Suppose $L_1,\ldots, L_m$ are the local times of $m\ge 2$ independent Brownian motions,
and $q_j\ge 1$ for all $1\le j \le m$. Then, for $q := \sum_{j=1}^m q_j$,
\begin{itemize}
\item [(a)] $\displaystyle \prob\Big\{ \int_{-\infty}^\infty \prod_{i=1}^m L_i^{q_i}(x,\sigma^{\ssup i})\, dx
< \eps \Big\} \asymp \eps^{\frac{2}{1+q}}\, ,$\\[1mm]
\item [(b)] $\displaystyle \prob\Big\{ \int_{-\infty}^\infty \prod_{i=1}^m L_i^{q_i}(x,1)\, dx
< \eps \Big\} \asymp  \eps^{\frac{2}{1+q}}\, .$\\
\end{itemize}
\end{thm}

{\bf Remark:} The excluded case $m=1$ is entirely different, as the small value probabilities decay exponentially.
This will be discussed in Section~\ref{self} using the technique of the B\"ottcher case.\qed
\smallskip
\pagebreak[2]

Before giving the detailed proof we show how the analogy to the martingale limit of
a Galton-Watson tree arises. From the Brownian paths we need to recognise the 
particular elements of the tree featuring in the proof of the Schr\"oder case: 
For each vertex of the spine we first need to decide whether a subtree splits off from 
the vertex (this happens independently with probability~$1-p_1$), and supposing 
this happens at the vertex in the $k^{\rm th}$ generation, we need to see that 
this subtree gives rise to a summand of the intersection local time,
which in distribution equals $\mu^{-k}$ times the intersection local time. 
Once an inequality analogous to \eqref{spinaldecomp} is
established
, we get lower tail asymptotics featuring
the parameters $\mu$ and $p_1$ used in the construction of the tree.

\begin{figure}[htbp]
\centering
\includegraphics[width=0.65\textwidth]{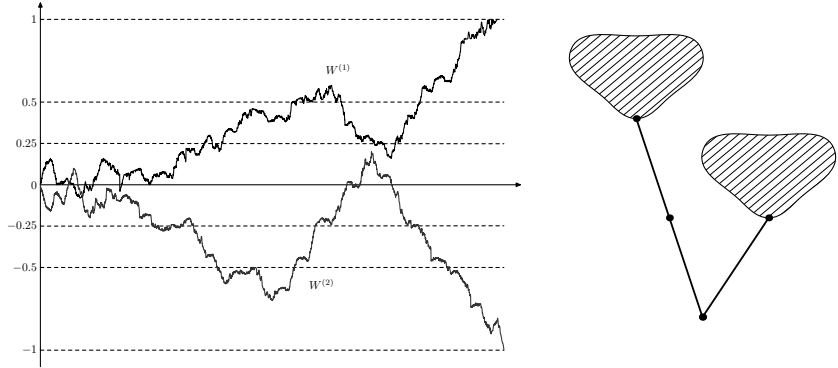} 
\label{tree}
\caption{The tree associated to two Brownian paths for $\eta=2$, up to $2^{\rm nd}$ generation. 
The intervals $W^{\ssup 1}[\tau^{\ssup 1}_{1},\tau^{\ssup 1}_{0}]$ and 
$W^{\ssup 2}[\tau^{\ssup 2}_{1}, \tau^{\ssup 2}_{0}]$ have a nonempty  intersection, 
and therefore the root has more than one offspring; by contrast the intervals $W^{\ssup 1}[\tau^{\ssup 1}_{2},\tau^{\ssup 1}_{1}]$ 
and $W^{\ssup 2}[\tau^{\ssup 2}_{2}, \tau^{\ssup 2}_{1}]$ are disjoint and therefore the second vertex on the spine has just one offspring.}
\end{figure}

To sketch the actual construction, focusing on~$m=2$ for the moment, we let $W^{\ssup 1}, W^{\ssup 2}$ 
be two independent Brownian motions started at the origin and assume that $W^{\ssup 1}$ exits $(-1,1)$ at the 
upper, and $W^{\ssup 2}$ exits $(-1,1)$ at the lower end of the interval. Fix 
$\eta>1$ and divide the Brownian paths according to the stopping~times
$$\tau^{\ssup 1}_k:=\inf\big\{ t\ge 0 \, : \, W^{\ssup 1}(t)= \eta^{-k} \big\}\qquad \mbox{ and } \qquad
\tau^{\ssup 2}_k:=\inf\big\{ t\ge 0 \, : \, W^{\ssup 2}(t)= -\eta^{-k} \big\}\, .$$
To build the tree from its spine $v_0(1),\ldots,v_n(1)$ of leftmost particles in the first $n$ generations,
we let the $k^{\rm th}$ individual $v_0(k)$ on this spine have more than one offspring if
$$W^{\ssup 1}[\tau^{\ssup 1}_{k+1},\tau^{\ssup 1}_{k}] \cap W^{\ssup 2}[\tau^{\ssup 2}_{k+1},\tau^{\ssup 2}_{k}]\not= \emptyset\, .$$
If the intervals intersect, the intersection local time of the two Brownian motions $W^{\ssup j}$, started at time
$\tau^{\ssup j}_{^{k+1}}$ and stopped at the time $\tau^{\ssup j}_{^k}$, for $j\in\{1,2\}$, give rise to a summand of the 
total intersection local time which is approximately distributed like a scaled copy of the total intersection local time.
\pagebreak[2]

\subsection{Intersection local times: The parameters $\mu$ and $p_1$}

We start with a basic scaling property of intersection local times. For any 
points $x_1,\ldots,x_m\in\R$ we suppose that under $\prob_{(x_j)}$ the
Brownian motion $W^{\ssup j}$ is started in $x_j$, and for $\eta>0$ we denote by
$$\tau^{\ssup j}(\eta)=\inf\{t>0 \colon W^{\ssup j}(t)=\eta\}$$
the first hitting time of $\eta$ by the Brownian motion $W^{\ssup j}$.

\begin{lemma}\label{scaling}
For every $\eps>0$ and for $q := \sum_{j=1}^m q_j$ we have
$$\prob_{(x_j/\eta)}\Big\{ 
\int_{-\infty}^\infty \prod_{j=1}^m L_j^{q_j}\big(x, \tau^{\ssup j}(1)\big) \, dx  <\eps  \Big\}
= \prob_{(x_j)}\Big\{  \int_{-\infty}^\infty \prod_{j=1}^m L_j^{q_j}\big(x, \tau^{\ssup j}(\eta)\big) \, dx 
<\eps  \eta^{1+q}\Big\}.$$
\end{lemma}

\begin{proof}
By Brownian scaling we have 
$$\prob_{x_j/\eta}\Big\{ L(x, \tau^{\ssup j}(1))<\eps \Big\}=
\prob_{x_j}\Big\{ \eta^{-1}\, L(\eta x, \tau^{\ssup j}(\eta))< \eps \Big\}.$$
Hence 
$$\begin{aligned}
\prob_{(x_j/\eta)}\Big\{  & \int_{-\infty}^\infty \prod_{j=1}^m L_j^{q_j}\big(x, \tau^{\ssup j}(1)\big) \, dx 
<\eps \Big\} =\prob_{(x_j)}\Big\{ \eta^{-\sum_{j=1}^m q_j}  
\int_{-\infty}^\infty \prod_{j=1}^m L_j^{q_j}\big(\eta x, \tau^{\ssup j}(\eta)\big) \, dx <\eps \Big\}\\
& = \prob_{(x_j)}\Big\{   \eta^{-(1+q)} \,\int_{-\infty}^\infty \prod_{j=1}^m L_j^{q_j}
\big(x, \tau^{\ssup j}(\eta)\big) \, dx <\eps\Big\},
\end{aligned}$$
and this proves the lemma.
\end{proof}

Fix $\eta>1$ and let $W^{\ssup 1}, \ldots, W^{\ssup m}$ be Brownian motions started in
the origin. Fix a set $M\subset\{1,\ldots, m\}$ and define stopping times
$$\tau^{\ssup j}_k:=\tau^{\ssup j}_k(M) := \left\{
\begin{array}{ll} \inf\big\{ t\ge 0 \, : \, W^{\ssup j}(t)=\eta^{-k} \}\,  &
\mbox{ if } j\in M\, ,\\
 \inf\big\{ t\ge 0 \, : \, W^{\ssup j}(t)=-\eta^{-k} \}\,  &
\mbox{ if } j\not\in M\, ,\\
\end{array}\right.$$
and abbreviate $\tau^{\ssup j}:=\tau^{\ssup j}_0(M)$. Suppose that under $\p_{(\pm \eps)}$ the Brownian motion 
$W^{\ssup j}$ is started in the point $+\eps$, if $j\in M$, and in the point $-\eps$ otherwise.
\pagebreak[3]

For $0<s<t$, define local times $L_j(x,s,t):=L_j(x,t) - L_j(x,s)$ over the time interval $[s,t]$, and
$$L_k := \int_{-\infty}^\infty \prod_{j=1}^m L_j^{q_j}\big(x,\tau^{\ssup j}_{k+1}, \tau^{\ssup j}_{k}\big) \, dx\, .$$
By the previous lemma, for every $k$, we have 
\begin{equation}\label{scalingvar}
\eta^{k(1+q)}\, L_k \stackrel{d}{=} L_0.
\end{equation}
This identifies the parameter $\mu$ 
as $\eta^{1+q}$. Recall that in the tree model this parameter corresponds to the mean offspring number.

\begin{lemma}\label{findp1}
If $M$ is a proper, nonempty subset of $\{1,\ldots,m\}$, we have
$$\prob_{(\pm \eps)}\big\{ W^{\ssup 1}[0,\tau^{\ssup 1}] \cap \ldots \cap W^{\ssup m}[0,\tau^{\ssup m}] = \emptyset \big\}
\asymp \eps^2\, .$$
\end{lemma}

\begin{proof}
On the one hand, if $\{ W^{\ssup 1}[0,\tau^{\ssup 1}] \cap \ldots \cap W^{\ssup m}[0,\tau^{\ssup m}]  = \emptyset \}$,
then at least one of the motions $W^{\ssup j}$, $j\in M$, does not reach level $-\eps$ before level $1$,
the probability of this being $2\eps/(1+\eps)$ per motion by the gambler's ruin probability. 
Analogously, one of the motions $W^{\ssup j}$, $j\not\in M$, does not reach level $\eps$ before level $-1$, which has the same probability. This gives the upper bound
$$\prob_{(\pm \eps)}\big\{ W^{\ssup 1}[0,\tau^{\ssup 1}] \cap \ldots \cap W^{\ssup m}[0,\tau^{\ssup m}]  
= \emptyset \big\} \le   \sfrac{\eps^2}{(1+\eps)^2}  \, 4\ell(m-\ell)\, ,$$
where $\ell$ is the cardinality of $M$. For the lower bound, note that if one of the motions in each of the two groups
does not reach level $0$ before level $1$, resp. $-1$, this implies
$W^{\ssup 1}[0,\tau^{\ssup 1}] \cap \ldots \cap W^{\ssup m}[0,\tau^{\ssup m}] = \emptyset$. As, for each motion,
this event has probability $\eps$, we obtain
$$\prob_{(\pm \eps)}\big\{  W^{\ssup 1}[0,\tau^{\ssup 1}] \cap \ldots \cap W^{\ssup m}[0,\tau^{\ssup m}] = \emptyset \big\} \ge \eps^2\, . \\ \ \vspace{-0.7cm}$$
\end{proof}

\smallskip

{\bf Remark:}
A refined calculation along the same lines shows that, as $\eps\downarrow 0$,
$$\prob_{(\pm \eps)}
\big\{W^{\ssup 1}[0,\tau^{\ssup 1}] \cap \ldots \cap W^{\ssup m}[0,\tau^{\ssup m}]  = \emptyset \big\}
\sim \eps^2\,\, 2\ell(m-\ell)  \, ,$$
where $\ell$ is the cardinality of $M$, but we do not need this here.\qed

\smallskip

By Brownian scaling we infer from Lemma~\ref{findp1} that there are constants $0<c<C$ such
that, if $M\subset\{1,\ldots,m\}$ is proper and nonempty, for any nonnegative integer $k$ and $\eta>1$,
$$c\,\eta^{-2} \le \prob\big\{ W^{\ssup 1}[\tau_{k+1}^{\ssup 1},\tau_k^{\ssup 1}] \cap \ldots \cap 
W^{\ssup m}[\tau_{k+1}^{\ssup 1},\tau_{k}^{\ssup m}] = \emptyset \big\} \le C \, \eta^{-2}\, ,$$
and thus the parameter $p_1$ is identified (with sufficient accuracy) as $\eta^{-2}$. Recall that $p_1$ corresponds
in the tree model to the probability that a vertex has only one offspring.

\subsection{Intersection local times: The lower bound}

Let $W^{\ssup 1}, \ldots, W^{\ssup m}$ be Brownian motions started at the origin, and
fix $M\subset\{1,\ldots,m\}$ such that $1\in M$ and $2\not\in M$. We propose a sufficient 
strategy to realise the event  $\{X(\sigma^{\ssup 1}(1),\ldots, \sigma^{\ssup m}(1)) <\eps\}$,
which is time-inhomogeneous and consists of two phases. Given $\eps>0$ the phases are 
separated by the stopping times
$$\omega^{\ssup j} := \inf\big\{ t\ge 0 \colon W^{\ssup j} \not\in(-\eps^{1/(1+q)}, \eps^{1/(1+q)}) \big\}, \quad \mbox{ for } j\in\{1,\ldots,m\} .$$
The first phase is described by the event
$$\begin{aligned}
E_1: = \big\{ W^{\ssup j}(\omega^{\ssup j})= \pm & \eps^{1/(1+q)}, \,
\inf\{ \pm W^{\ssup j}(s) \colon 0\le s\le\omega^{(j)} \} > -\sfrac 12 \eps^{1/(1+q)}\\
&\mbox{ for all $j$ and } X(\omega^{\ssup 1},\ldots, \omega^{\ssup m}) <\eps \big\},
\end{aligned}$$
where $\pm$ indicates $+$ if $j\in M$ and $-$ otherwise.
By the scaling verified in Lemma~\ref{scaling} the probability $\delta:=\prob(E_1)>0$ does not
depend on $\eps$. The second phase is described by the event
$$\begin{aligned}
E_2: = \big\{ W^{\ssup j}(\tau^{\ssup j}) & = \pm 1\mbox{ for all $j$ and }
\inf\{ W^{\ssup 1}(s) \colon \omega^{\ssup 1} \le s\le \tau^{\ssup 1} \} \ge \sfrac 12 \,
\eps^{1/(1+q)}, \\  & \mbox{ and }
\sup\{ W^{\ssup 2}(s) \colon \omega^{\ssup 2} \le s\le \tau^{\ssup 2} \}
\le  -\sfrac 12\, \eps^{1/(1+q)} \big\} .
\end{aligned}$$
Observe that, if $E_1$ and $E_2$ hold, we have
$$X(\sigma^{\ssup 1},\ldots, \sigma^{\ssup m})=X(\tau^{\ssup 1},\ldots, \tau^{\ssup m})=
X(\omega^{\ssup 1},\ldots, \omega^{\ssup m}) <\eps\, ,$$
as required. Moreover, using the strong Markov property and the gambler's ruin estimate,
$$\begin{aligned}
\prob( E_1 \cap E_2) & = \me\big[ \1_{E_1} \, \prob_{{(W^{\ssup j}(\omega^{\ssup j}))}}( E_2) \big] \\[1mm]
& = \prob(E_1) \, 
\bigg( \frac{1+ \eps^{1/(1+q)}}{2}\bigg)^{m-2}
\bigg( \frac{\frac 12 \eps^{1/(1+q)}}{1-\frac 12 \eps^{1/(1+q)}}\bigg)^2
\, ,
\end{aligned}$$
so the lower bound holds with $c:=\delta (1/2)^{m}$.

\subsection{Intersection local times: The logarithmic upper bound}

We now give an upper bound for the small value probability
of $X(\sigma^{\ssup 1},\ldots, \sigma^{\ssup m})$ along the lines of the
argument leading to \eqref{roughbound}. Fix an arbitrarily
small $\delta>0$. Let $C\ge 1$ be the constant in the implied upper
bound of Lemma~\ref{findp1}. Choose and fix an integer $\eta>(2C)^{1/\delta}$.

For any subset $M\subset \{1,\ldots, m\}$ define the event
$$E(M):= \big\{ W^{\ssup j}(\sigma^{\ssup j})=1 \mbox{ for all }j\in M, \,
W^{\ssup j}(\sigma^{\ssup j})=-1 \mbox{ for all }j\not\in M \big\}\, .$$
Recall the definition of the stopping times $\tau^{\ssup j}_k := \tau^{\ssup j}_k(M)$.
Then
\begin{equation}\begin{aligned}
\prob\big\{ & X(\sigma^{\ssup 1},\ldots, \sigma^{\ssup m}) < \eps \big\}
& = \sum_{M\subset\{1,\dots,m\}}
\prob\Big( \{ & X(\tau^{\ssup 1}_0,\ldots, \tau^{\ssup m}_0)< \eps \} \cap E(M) \, \Big)\, .
\label{stopping}
\end{aligned}\end{equation}
It therefore suffices to fix $M\subset \{1,\ldots, m\}$ and give upper bounds
for $\prob\{X(\tau^{\ssup 1}_0,\ldots, \tau^{\ssup m}_0)< \eps \}$. Define, for $0<s<t$, local times
$L_j(x,s,t):=L_j(x,t) - L_j(x,s)$ over the time interval $[s,t]$. Denote
$$L_k := \int_{-\infty}^\infty \prod_{j=1}^m L_j^{q_j}(x,\tau^{\ssup j}_{k+1}, \tau^{\ssup j}_{k}) \, dx\, .$$
Then the random variables $X_k= \eta^{k(1+q)}\, L_k$  are independent,
by the Markov property, and identically distributed, by~\eqref{scalingvar}.
By Lemma~\ref{findp1} we have $\prob\{ X_0=0 \} \le C \eta^{-2}$ if $M$ is a proper, nonempty
subset of $\{1,\ldots,m\}$, and otherwise obviously $\prob\{ X_0=0 \}=0$. This implies that there
exists a $\theta>0$ such that
$$\prob\{ X_0< \theta \} \le 2C\, \eta^{-2}\, .$$
Now, given $\eps>0$ pick the integer~$n$ such that
$$\theta\, \eta^{-(n+1) (1+ q)}< \eps\le \theta\, \eta^{-n (1+ q)}\, ,$$
Note that for $q_i \geq 1$, by super-additivity of $x\mapsto x^{q_i}$, $x\ge 0$ we get
\[ L_j^{q_j} (x,\tau^{\ssup j}_0 ) \geq 
\bigg( \sum_{k=0}^{n-1} L_j (x, \tau_{k+1}^{\ssup j}, \tau_k^{\ssup j}) \bigg)^{q_j} \geq \sum_{k=0}^{n-1} L_j^{q_j} (x, \tau_{k+1}^{\ssup j}, \tau_k^{\ssup j}) .  \]
Applying this to the intersection local times, it follows that,
\[  X(\tau^{\ssup 1}_0, \ldots\ ,\tau^{\ssup m}_0) = \int_{-\infty}^\infty \prod_{j=1}^m L_j^{q_j} (x, \tau_0^{\ssup j}) \, dx
  \geq  \int_{-\infty}^\infty \prod_{j=1}^m \bigg( \sum_{k=0}^{n-1} L_j^{q_j} (x, \tau_{k+1}^{\ssup j}, \tau_k^{\ssup j})\bigg) \, dx \geq
\sum_{k=0}^{n-1} L_k . \]
Hence we can estimate
$$\begin{aligned}
\prob\big\{  X(\tau^{\ssup 1}_0, \ldots\ , \tau^{\ssup m}_0) < \eps \big\} &\leq \prob\Big\{\sum_{k=0}^{n-1} L_k < \eps \Big\}
\leq \prob\Big\{\sum_{k=0}^{n-1} \eta^{-k (1+ q)}X_k< \theta\, \eta^{-n (1+ q)}\Big\}\\
&\leq  \prob\Big\{\sum_{k=0}^{n-1} X_k < \theta \Big\} \leq \big(\prob\{X_0<\theta\}\big)^{n}
\leq \, \big(2C\big)^{n}  \eta^{-2n} \le K\, \eps^{\frac{2-\delta}{1+ q}}\, ,
\end{aligned}$$
for the constant $K:=\eta^{2-\delta} \theta^{\frac{-2+\delta}{1+ q}}$.
As $\delta>0$ can be chosen arbitrarily small, this shows that
\begin{equation}\label{rough2}
\limsup_{\eps\downarrow 0} \frac{\log\prob\big\{ X(\sigma^{\ssup 1},\ldots, \sigma^{\ssup m}) < \eps \big\}}{-\log \eps}
\le {\frac{-2}{1+ q}}\, .\end{equation}
Note (for use in Lemma~\ref{beta}) that the proof also shows that \eqref{rough2} holds if 
$W^{\ssup 1},\ldots, W^{\ssup m}$ are started in arbitrary points of the interval
$[-\eta^{-n},\eta^n]$ instead of the origin.

\subsection{Intersection local times: Up-to-constant asymptotics}

Fix the set $M\subset\{1,\ldots,m\}$, the integer $\eta>1$, and recall the notation from the previous
section. Define a sequence $(a(n) \colon n\ge 0)$ by 
$$a(n) := \p\{X(\tau^{\ssup 1}_0,\ldots,\tau^{\ssup m}_0) < \theta\eta^{-n(1+q)}\}\, \eta^{2n}\, .$$
Given $0<\eps<1$ we again pick the integer $n$ such that $\theta\eta^{-(n+1)(1+q)}\le \eps <
\theta\eta^{-n(1+q)}$. Then
$$\begin{aligned}
\prob\big\{ X(\tau^{\ssup 1}_0,\ldots,\tau^{\ssup m}_0) < \eps \big\} 
& \le \p\big\{ X(\tau^{\ssup 1}_0,\ldots,\tau^{\ssup m}_0) <
\theta\eta^{-n(1+q)} \big\} \\
& = a(n)\,\eta^{-2n} \le a(n) \, \eta^2\, \theta^{-\frac{2}{1+q}}\, \eps^{\frac{2}{1+ q}}\, , 
\end{aligned}$$
hence, to complete the proof, it suffices to show that $(a(n) \colon n\ge 0)$ is bounded.
Define $$T:=\min\big\{ k\ge 0\colon W^{\ssup 1}[\tau^{\ssup 1}_{k+1},\tau^{\ssup 1}_{0}]\cap \ldots\cap
W^{\ssup m}[\tau^{\ssup m}_{k+1},\tau^{\ssup m}_{0}] \not=\emptyset \big\}\, .$$
In our tree heuristic $T$ is the first generation in which a tree is branching off the spine.
The next lemma controls the behaviour of this tree and plays a similar r\^ole to \eqref{estimate_for_W_fixed_T}. 

\begin{lemma}\label{beta}
There exists a sequence $(\beta(i)\colon i\in \N)$ of nonnegative numbers
with $\sum \beta(i)<\infty$ such that, for $0\le j\le n-1$,
$$\prob_{(y_i)} \big\{ X(\tau^{\ssup 1}_0,\ldots,\tau^{\ssup m}_0) < \theta\eta^{-n(1+q)},
\, T=j  \big\} \le \eta^{-2j-2}\, \beta(n-j-1)\, ,$$
where $y_i=\pm\,\eta^{-j-1}$ with the sign chosen according to whether $i\in M$ or not.
\end{lemma}

\begin{proof}
For $i\in\{1,\ldots,m\}$ and $k\in\{-\eta^{n-j-1}, \ldots, \eta^{n-j-1}-1\}$,
we introduce stopping times, 
$$\varrho^{\ssup i}_k:=\inf\big\{ t\ge 0 \, : \, W^{\ssup i}(t)\in[k \eta^{-n}, (k+1)\eta^{-n}] \big\}\, .$$
The assumption $T=j$ implies that there exists $k\in\{-\eta^{n-j-1}, \ldots, \eta^{n-j-1}-1\}$
such that $\varrho^{\ssup i}_k<\tau^{\ssup i}_0$, for all $i\in\{1,\ldots, m\}$.
If this holds, then let $\sigma^{\ssup i}_{j}:=\inf\{ t\ge\varrho^{\ssup i}_k
\, : \, W^{\ssup i}(t)=\pm \eta^{-j}\}$ (with the usual convention on $\pm$). Hence, for
any $0<\delta<1$ and sufficiently large $n-j$, using first Lemma~\ref{findp1}
with $\eps=\eta^{-j}$, then \eqref{rough2} and the subsequent remark in combination with
Lemma~\ref{scaling} and, of course, the strong Markov property,
$$\begin{aligned}
\prob_{(y_i)} & \big\{ X(\tau^{\ssup 1}_0,\ldots,\tau^{\ssup m}_0) < \theta\eta^{-n(1+q)},
\, T=j  \big\} \\
& \le  \sum_{k=-\eta^{n-j-1}}^{\eta^{n-j-1}-1}
\me_{(y_i)}\Big[ \1\big\{X(\sigma^{\ssup 1}_{j},\ldots,\sigma^{\ssup m}_{j}) < \theta\eta^{-n(1+q)}\big\}\\
& \qquad\qquad\qquad \times \, \prob_{(W^{\ssup i}(\sigma^{\ssup i}_{j}))}\big\{W^{\ssup 1}[0,\tau^{\ssup 1}_{0}]
\cap \ldots\cap W^{\ssup m}[0,\tau^{\ssup m}_{0}] =\emptyset \big\}  \Big] \\
& \le  \sum_{k=-\eta^{n-j-1}}^{\eta^{n-j-1}-1}
\prob_{(W^{\ssup i}(\varrho^{\ssup i}_{k}))}\big\{X(\sigma^{\ssup 1}_{j},\ldots,\sigma^{\ssup m}_{j}) 
< \theta\eta^{-n(1+q)} \big\} \; C\,\eta^{-2j}\\
& \le   2\eta^{n-j-1} \;\; \eta^{(-2+\delta)(n-j)} \; C\,\eta^{-2j}\, ,
\end{aligned}$$
which gives the result with  $\beta(i):=2C\,\eta^{\delta}\, \eta^{(-1+\delta)i}$.
\end{proof}

We now argue as in \eqref{inequality_for_W} of the Schr\"oder case, using in the second step the upper
bound of Lemma~\ref{findp1} and denoting the implied constant there by $C>0$,
\begin{equation}\begin{aligned}\label{start}
\prob& \big\{ X(\tau^{\ssup 1}_0,\ldots,\tau^{\ssup m}_0) < \theta\eta^{-n(1+q)} \big\} \\
& \le \prob\big\{ T \ge n \big\}  + \sum_{j=0}^{n-1}
\prob\big\{ X(\tau^{\ssup 1}_0,\ldots,\tau^{\ssup m}_0) < \theta\eta^{-n(1+q)}, \, T=j  \big\} \\
& \le C\, \eta^{-2n} + \sum_{j=0}^{n-1}
\prob\big\{ X(\tau^{\ssup 1}_0,\ldots,\tau^{\ssup m}_0) < \theta\eta^{-n(1+q)}, \, T=j  \big\} \, .
\end{aligned}\end{equation}
To estimate the remaining probability we first use the strong Markov property,
then Lemma~\ref{beta} to estimate the inner probability, and finally 
the definition of $(a(n) \colon n\ge0)$ in combination with Lemma~\ref{scaling}, to obtain
$$\begin{aligned}
\prob\big\{ X(\tau^{\ssup 1}_0,\ldots, & \tau^{\ssup m}_0) < \theta\eta^{-n(1+q)}, \, T=j  \big\} \\[1mm]
& \le  \me\Big\{ \1\{ X(\tau^{\ssup 1}_{j+1},\ldots,\tau^{\ssup m}_{j+1})< \theta\eta^{-n(1+q)}  \} \\
& \qquad\qquad\qquad \times\prob_{(W^{\ssup i}(\tau^{\ssup i}_{j+1}))} \big\{ X(\tau^{\ssup 1}_0,\ldots,\tau^{\ssup m}_0) < 
\theta\eta^{-n(1+q)},\, T=j  \big\} \Big\} \\[1mm]
& \le \eta^{-2j-2}\, \beta(n-j-1) \, \prob\big\{ X(\tau^{\ssup 1}_{j+1},\ldots,\tau^{\ssup m}_{j+1}) <
\theta\eta^{-n(1+q)} \big\} \\[1mm]
& \le \eta^{-2n}\, \beta(n-j-1) \, a(n-j-1) \, .\end{aligned}$$
Plugging this into \eqref{start} we obtain a recursion formula for $a(n)$, namely
$$a(n) \le \sum_{j=0}^{n-1} \beta(n-j-1) \, a(n-j-1)+C \qquad\mbox{ for } n\ge 0.$$
As before, boundedness of $(a(n) \colon n\ge 0)$ follows from the recursion and the
fact that $\sum \beta(j)<\infty$.

\subsection{Intersection local times at fixed times}

In this section we use a technique adapted from~\cite{La96} to transfer our results
from hitting times to fixed times, thus proving Theorem~\ref{main2}\,(b). Recall the
following simple tail estimates for the first exit times~$\sigma^{\ssup j}(x)$ from 
the interval $(-x,x)$ by a Brownian motion $W^{\ssup j}$ started in $x_j$. 
 
\begin{lemma}\label{min_max_estimate} 
There exist constants $\beta>0$ and $\kappa>0$ such that, for all $x > 0$, $|x_j| \leq x/2$ and $a>0$,
\begin{itemize}
\item[(a)] $\displaystyle \p_{(x_j)} \Big\{ \min_{j=1}^m \,\sigma^{\ssup j}(x) \leq a x^2 \Big\} \leq 
\kappa\, e^{- \beta/ a},$\\[2mm]
\item[(b)] $\displaystyle \p_{(x_j)} \Big\{ \max_{j=1}^m \,\sigma^{\ssup j}(x) \geq a x^2 \Big\} \leq 
\kappa\, e^{- \beta a} \, . $
\end{itemize}
\end{lemma}

\begin{proof}
By scaling, we may assume that~$x =1$. On the one hand, using the reflection principle, we~get
$$\p_{x_j} \{ \sigma^{\ssup j}(1) \leq a\} \leq \p_0 \big\{ \sup_{t \leq a} |W^{\ssup j}(t)| \geq \tfrac{1}{2} 
\big\} \leq 2\,  \p_{0} \{ |W^{\ssup j}(a)| \geq \tfrac{1}{2} \} = 2 \p_{0} \big\{ |W^{\ssup j}(1)| 
\geq \tfrac{1}{2\sqrt{a}} \big\} \, ,$$
and hence (a) follows from a standard estimate for the tail of a normal distribution.
On the other hand, (b) follows from 
$\p_{x_j}\{ \sigma^{\ssup j}(1) \geq k \, \big| \, \sigma^{\ssup j}(1) \geq k-1\} 
\leq \p_{0} \{ |W^{\ssup j}(1)| \leq 2 \} < 1$ by iteration. 
\end{proof}

For the \emph{lower bound} we get, for any $a>0$, using Lemma~\ref{scaling} in the second step,
$$\begin{aligned} 
\p & \{ X(1,\ldots,1) < \eps \}  \\ & 
\geq \p\big\{ X(\sigma^{\ssup 1}(a), \ldots, \sigma^{\ssup m}(a)) < \eps \big\} - 
\p\big\{ X(\sigma^{\ssup 1}(a), \ldots, \sigma^{\ssup m}(a))<\eps, \, \min_{j=1}^m \sigma^{\ssup j}(a) \leq 1 \big\} \\
& = \p\big\{ X(\sigma^{\ssup 1}(1), \ldots, \sigma^{\ssup m}(1)) < a^{-(1 + q)} \eps \big\} 
- \p\big\{ X(\sigma^{\ssup 1}(a), \ldots, \sigma^{\ssup m}(a)) < \eps, \, \min_{j=1}^m 
\sigma^{\ssup j}(a) \leq 1 \big\} \, . \end{aligned} $$
Using first Theorem~\ref{main2}\,(a) in combination with Lemma~\ref{scaling}
and then Lemma~\ref{min_max_estimate}\,(a), 
\[\begin{aligned} \p \big\{ & X(\sigma^{\ssup 1}(a), \ldots, \sigma^{\ssup m}(a)) < \eps, \, 
\min_{j=1}^m \sigma^{\ssup j}(a) \leq 1 \big\} \\
 & \leq \me \Big[ \1\big\{ X(\sigma^{\ssup 1}(a/2), \ldots, \sigma^{\ssup m}(a/2)) < \eps \big\}  \, 
\p_{(W^{(j)}(\sigma^{\ssup j}(a/2)))} \big\{  \min_{j=1}^m \sigma^{\ssup j}(a) \leq 1 \big\} \Big] \\
 & \leq 4C a^{-2} \eps^{\frac{2}{1 + q}} \, \sup_{|x_j| = a/2} \p_{(x_j)} 
\big \{  \min_{j=1}^m \sigma^{\ssup j}(a)  \leq 1 \big\} 
\leq 4C a^{-2} \eps^{\frac{2}{1 + q}} \kappa\, e^{- \beta a^2} \, , \end{aligned} \]
where $C>0$ is the implied constant in the upper bound of Theorem~\ref{main2}\,(a). Substituting this 
into the previous equation and applying the lower bound of Theorem \ref{main2}\,(a) 
with the implied constant denoted by $c>0$, we get 
\[ \begin{aligned} \p \{  X(1,\ldots,1)  < \eps \} 
& \geq \p \{ X(\sigma^{\ssup 1}(1), \ldots, \sigma^{\ssup m}(1)) < a^{-(1 + q)} \eps \} 
- 4C a^{-2} \eps^{\frac{2}{1 + q}} \kappa\, e^{- \beta a^2} \\
& \geq \big( c a^{-2} - 4C a^{-2} \kappa e^{-\beta a^2} \big)\, \eps^{\frac{2}{1 + q}}\, , \end{aligned} \]
and the result follows if we choose $a$ large enough to ensure that the bracket is positive.


For the \emph{upper bound}, given $\eps>0$, we pick the integer $n$ such that
\begin{equation}\label{chon}
e^{- \beta 2^n} \leq \eps^{\frac{2}{1+q}} < e^{- \beta{2^{n-1}}}.
\end{equation}
We base the argument on the decomposition
\begin{equation} \begin{aligned} 
\p\{  X(1,\ldots,1) < \eps \} \leq \p \{ & X(\sigma^{\ssup 1}(1), \ldots, \sigma^{\ssup m}(1)) < \eps \} \\
& + \sum_{i=0}^{n-1} \sum_{j=1}^m \p\big\{ X(\sigma^{\ssup 1}(2^{-i-1}), \ldots, \sigma^{\ssup m}(2^{-i-1})) < \eps  \, , 
\sigma^{\ssup j}(2^{-i}) \geq 1 \big\} \\
& + \p \big\{ \max_{j=1}^m \sigma^{\ssup j} (2^{-n})\geq 1 \big\}. 
\end{aligned} \label{eqn:upper_bound_fixed_times}\end{equation}
We bound the first term on the right hand side using Theorem~\ref{main2}\,(a) and the last one 
using Lemma~\ref{min_max_estimate}\,(b) and \eqref{chon}. It remains to bound the sum in the middle. 
To this end we write
$$\sigma^{\ssup j}(2^{-i}) = \sum_{k=i}^{n} \big( \sigma^{\ssup j}(2^{-k}) - \sigma^{\ssup j}(2^{-(k+1)}) \big) 
+ \sigma^{\ssup j}(2^{-(n+1)}) \, , $$
and note that, as $2^{-2n-2}2^{n+i} + \sum_{k=i}^n2^{i-k-1}\le 1$, we get
$$\begin{aligned}
\p\big\{ & X(\sigma^{\ssup 1}(2^{-i-1}), \ldots, \sigma^{\ssup m}(2^{-i-1})) < \eps  \, , 
\sigma^{\ssup j}(2^{-i}) \geq 1 \big\} \\
& \le \sum_{k=i}^{n} \prob\{ X(\sigma^{\ssup 1}(2^{-i-1}), \ldots, \sigma^{\ssup m}(2^{-i-1})) < \eps  \, , 
 \sigma^{\ssup j}(2^{-k}) - \sigma^{\ssup j}(2^{-(k+1)}) \geq 2^{i-k-1} \big\}\\
& \qquad\qquad+ \prob\big\{ \sigma^{\ssup j}(2^{-(n+1)}) \geq 2^{-2n-2}2^{n+i} \big\}.
\end{aligned}$$
Again the contribution from the last summand can be bounded using Lemma~\ref{min_max_estimate}\,(b). For the
remaining term we use the strong Markov property to obtain, if $n \geq k \geq i+1$, 
\begin{equation} \begin{aligned} \p \{ X(& \sigma^{\ssup 1}(2^{-i-1}),   
\ldots, \sigma^{\ssup m}(2^{-i-1})) < \eps , \, \sigma^{\ssup j}(2^{-k}) - \sigma^{\ssup j}(2^{-k-1}) 
\geq 2^{i-k-1} \} \\[1mm]
& \leq \p \{ X(\sigma^{\ssup 1}(2^{-k-1}),  \ldots, \sigma^{\ssup m}(2^{-k-1})) < \eps \} 
\sup_{|x_j|= 2^{-k-1}} \p_{x_j}\{ \sigma^{\ssup j}(2^{-k}) \geq 2^{i-k-1} \} \\ 
& \quad \times 
\sup_{|x_j| = 2^{-k}} \p_{(x_j)} \{ X(\sigma^{\ssup 1}(2^{-i-1}),  \ldots, \sigma^{\ssup m}(2^{-i-1})) < \eps \}. 
\end{aligned} \label{eqn:split_up_interval} \end{equation}
If $n$ is large enough (or, equivalently, $\eps>0$ small enough) to satisfy $e^{-\beta 2^{n-2}}\le 2^{-n}$,
then we get that
$$\begin{aligned}
\sup_{|x_j| = 2^{-k}} & \p_{(x_j)} \big\{ X(\sigma^{\ssup 1}(2^{-i-1}),  \ldots, \sigma^{\ssup m}(2^{-i-1})) < \eps \big\}\\[-1mm]
& = \sup_{|x_j| = 2^{i-k+1}} \p_{(x_j)} \big\{ X(\sigma^{\ssup 1}(1),  \ldots, \sigma^{\ssup m}(1)) < 
\eps \, 2^{(i+1)(1 + q)} \big\}\\
& \le \sup_{|x_j| = 2^{i-k+1}} \p_{(x_j)} \big\{ X(\sigma^{\ssup 1}(1),  \ldots, \sigma^{\ssup m}(1)) < 
2^{(i-k+1)(1 + q)} \big\}\, .
\end{aligned}$$
Recall that $\tau^{\ssup j}(x)=\inf\{ t\ge 0 \colon W^{\ssup j}(t)=x\}$ and note that, for $|x_j|= 2^{-k}$,
$$\begin{aligned}
\p_{(x_j)} & \big\{ X(\sigma^{\ssup 1}(1),  \ldots, \sigma^{\ssup m}(1)) < 2^{(i-k+1)(1 + q)} \big\} \\
& \le \prob_{(2^{i-k+1})} \big\{ X(\sigma^{\ssup 1}(1),  \ldots, \sigma^{\ssup m}(1)) < 2^{(i-k+1)(1 + q)} \big\}\\
& \qquad + \prob_{(-2^{i-k+1})} \big\{ X(\sigma^{\ssup 1}(1),  \ldots, \sigma^{\ssup m}(1)) < 2^{(i-k+1)(1 + q)} \big\}\\
&\qquad + \sum_{j=1}^m\sum_{\ell=1}^m \p_{x_j}\big\{ \tau^{\ssup j}(2^{i-k+1})>\sigma^{\ssup j}(1) \big\}
\p_{x_\ell}\big\{ \tau^{\ssup \ell}(-2^{i-k+1})>\sigma^{\ssup \ell}(1) \big\}\, .
\end{aligned}$$
While the first two probabilities are bounded by constant multiples of $2^{2(i-k+1)}$
by Theorem~\ref{main2}\,(a), the double sum is bounded by $m^2\, 2^{2(i-k+2)}$ by the gambler's ruin probability.
Hence, for a suitable constant $C_0>1$ and all $n \geq k\ge i+1$,
$$\sup_{|x_j| = 2^{-k}} \p_{(x_j)} \big\{ X(\sigma^{\ssup 1}(2^{-i-1}),  \ldots, \sigma^{\ssup m}(2^{-i-1})) < \eps \big\}
\le C_0\, 2^{2(i-k)}\, .$$
Combining this with Lemma~\ref{min_max_estimate}\,(b) and substituting into~\eqref{eqn:split_up_interval} we get
for all $n \ge k\ge i$,
$$\begin{aligned} \p \big\{ X(& \sigma^{\ssup 1}(2^{-i-1}),   \ldots, \sigma^{\ssup m}(2^{-i-1})) < \eps , 
\, \sigma^{\ssup j}(2^{-k}) - \sigma^{\ssup j}(2^{-k-1}) \geq 2^{i-k-1} \big\} \\
& \leq C_1\, \eps^{\frac{2}{1 + q}} \big[ 2^{2k+2} e^{-\beta 2^{k+i-1}} 2^{2(i-k)} \big] 
\, , \end{aligned}$$ for $C_1:=C_0C\kappa$. After summing over $k\ge i$, $0\le i\le n-1$ and $1\le j \le m$, 
the square bracket on the right remains bounded, and this completes the proof of Theorem~\ref{main2}\,(b).

\section{Small value probabilities for self-intersection local times}\label{self}

In this section we look at a single Brownian motion and its $q$-fold self-intersection local time
$$X(t):=\int_{-\infty}^\infty L^{q}(x,t)\, dx\, .$$
This corresponds to the case $m=1$ of the scenario described in Section~\ref{mutual} and, as mentioned there,
this is quite different from the case $m>1$. The argument used to study the B\"ottcher case
of the Galton-Watson limit can be used to give an extremely simple proof of the following result.
\smallskip

\begin{thm}\label{main3}
Suppose $(L(x,t) \colon x\in\R,\, t\ge 0)$ is the local time field
and $\sigma:=\inf\{t\ge 0 \, : \, |B(t)|=1\}$ the first hitting time
of level one of a Brownian motion. Then, for every $q \geq 1$, we have
$$-\log \prob\Big\{ \int_{-\infty}^\infty L^{q}(x,\sigma)\, dx
< \eps \Big\} \asymp  \eps^{-\frac{1}{q}}\, .$$
\end{thm}

{\bf Remark:} The behaviour is radically different, when the Brownian motion is stopped
at a fixed time instead of a fixed level. Indeed, we will see in the proof of Theorem~\ref{main3}
that the optimal strategy to make $X(\sigma)$ small is simply to make $\sigma$ small, an option which cannot
be used to make $X(1)$ small. It was shown, for $q=2$ in \citet[Proposition~1]{HHK97} 
and extended to general $q>1$ by Xia Chen and Wenbo Li (unpublished), that there is a  constant $c(q)>0$ such that,
$$-\log\prob\Big\{ \int_{-\infty}^\infty L^q(x,1)\, dx < \eps \Big\}
\sim c(q)\, \eps^{\frac{-2}{q+1}}\, .\vspace{-8mm}$$\qed\\
\smallskip

\subsection{Self-intersection local time: The branching tree heuristic}

We first show how to establish the analogy between the
$q$-fold self-intersection local times and the martingale limit of a Galton-Watson
tree in the B\"ottcher case. The idea is to construct a nested family of random walks 
embedded into the Brownian path: The natural nesting of the embedded walks establishes 
the tree structure, and a constant multiple of the total number of steps of the finest embedded walk approximates 
the $q$-fold self-intersection local times. 
\smallskip

Let $(W(t) \colon t\ge 0)$ be a Brownian motion started at the origin and, for each nonnegative integer $n$,~let 
$$\mathfrak D_n := \big\{ k2^{-n} \, : \, k\in\{-2^n,\ldots, 2^n\}\big\} $$ 
be the collection of dyadic points of the $n^{\rm th}$ stage and let
$0=\tau_{^0}^{\ssup n}<\tau_{^1}^{\ssup n}<\cdots<\tau^{\ssup n}_{^{N(n)}}=\sigma$ 
be the collection of stopping times defined for $j\ge 1$ by
$$\tau_j^{\ssup n}:=\inf\big\{ t>\tau_{j-1}^{\ssup n} \, : \, W(t)\in\mathfrak{D}_n,\ \, 
W(t)\not= W(\tau_{j-1}^{\ssup n}) \big\}\, .$$
Then $(X^{\ssup n}(j) \colon 0 \le j \le N(n))$ defined by $$X^{\ssup n}(j):=2^n\, W(\tau_{j}^{\ssup n})$$
is the $n^{\rm th}$ embedded random walk and $N(n)$ its length.
We assign $N(1)$ offspring to the root, so that the vertices in the first generation 
correspond to the steps of height $1/2$ the path takes to reach level~$1$
or $-1$ for the first time. Then the number of children of each vertex in the first generation 
is determined by the number of steps of height $1/4$ the path makes during the step
of height $1/2$ corresponding to that vertex. This will be iterated ad infinitum to map
the Brownian path to an infinite tree. Note that the resulting tree is a Galton-Watson tree
and every vertex in this tree has at least two offspring, so that we are in the B\"ottcher case.

\begin{figure}[htbp]
\centering
\includegraphics[width=0.7\textwidth]{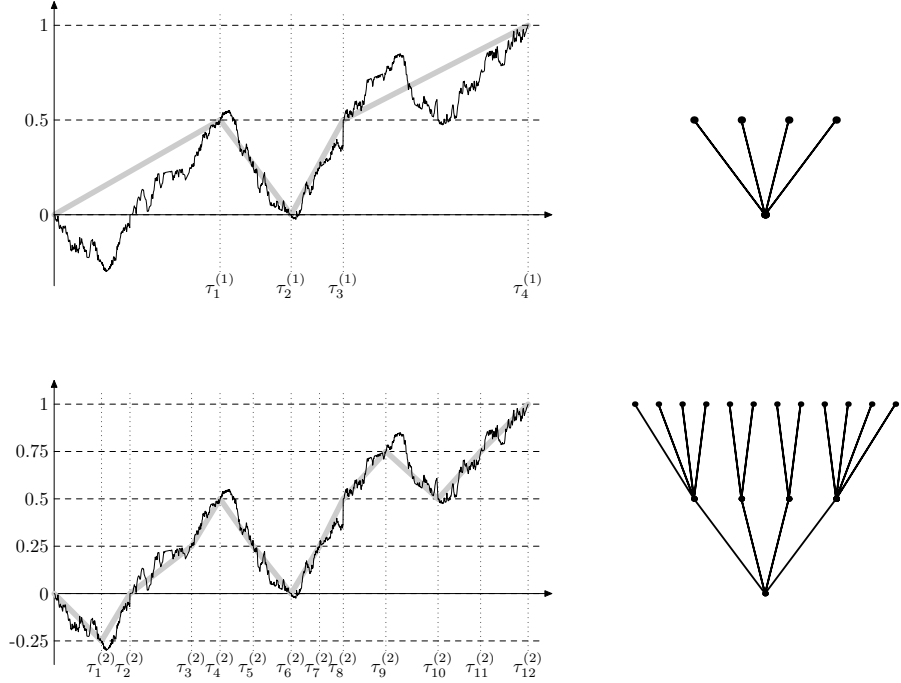} 
\caption{On the left, the first two embedded random walks with step sizes $\frac 12$, resp.~$\frac 14$,
on the right the corresponding first two generations of the associated tree.}
\end{figure}

\subsection{Self-intersection local time: The lower bound}

Recall from the last subsection the definition of the stopping times $0=\tau_0^{\ssup n}<\tau_1^{\ssup n}<\cdots
<\tau^{\ssup n}_{^{N(n)}}=\sigma$ and of $N(n)$. Note that $N(n)\ge 2^n$ and that $\prob\{ N(n)= 2^n\} =2\,(1/2)^{2^n}$.
Hence, for any $n$ and $\eps>0$,
$$\begin{aligned}
\prob\Big\{ \int_{-\infty}^\infty L^{q}(x,\sigma)\, dx
\le \eps \Big\} & \ge
\prob\Big\{ \int_{-\infty}^\infty L^{q}(x,\sigma)\, dx
\le \eps \,\Big| \, N(n) = 2^n\Big\}\times 2\,(1/2)^{2^n}\, .
\end{aligned}$$
By scaling, there exists a positive constant $C(q)$ such that, for
all $j\in\{1,\ldots, N(n)\}$, the random variables
$$Y_j:= C(q)\, 2^{n(1+q)}\int_{-\infty}^\infty L^{q}(x,\tau_{j-1}^{\ssup n}, \tau_{j}^{\ssup n})\, dx\, ,$$
have mean one. Given $\eps>0$ we pick the integer $n$ such that $2^{-(n+1)q}\le C(q) 2^{-2q}
\eps < 2^{-nq}$. Conditional on $N(n)=2^n$, for every $x\in\R$, we know that in the decomposition
$$L(x,\sigma) = \sum_{j=1}^{2^n} L(x,\tau_{j-1}^{\ssup n}, \tau_{j}^{\ssup n}) \, $$
only two summands can be non-zero. 
Thus, using the convexity of $x \mapsto x^q$ for $q \geq 1$, we obtain
$$ \int_{-\infty}^\infty L^{q}(x,\sigma)\, dx \le 2^{q-1} \, \sum_{j=1}^{2^n}
\int_{-\infty}^\infty L^{q}(x,\tau_{j-1}^{\ssup n}, \tau_{j}^{\ssup n})\, dx \le \, \eps
\, 2^{-1-n}\,\sum_{j=1}^{2^n} Y_j ,$$ and the summands on the right
are independent, identically distributed random variables with mean
one. Hence, by the law of large numbers,
$$\begin{aligned}
\prob\Big\{ \int_{-\infty}^\infty L^{q}(x,\sigma)\, dx
\le \eps \,\Big| \, N(n) = 2^n\Big\} & \ge
\prob\Big\{ 2^{-n}\,\sum_{j=1}^{2^n} Y_j
\le 2 \,\Big| \, N(n) = 2^n\Big\}
\stackrel{n\uparrow\infty}{\longrightarrow} 1\, ,
\end{aligned}$$
and, altogether, for $c(q):=4 \,(\log 2)\, C(q)^{-1/q} >0$ and all large values of $n$,
$$\prob\Big\{ \int_{-\infty}^\infty L^{q}(x,\sigma)\, dx
\le \eps \Big\} \ge (1/2)^{2^n} \ge \exp\big( -c(q) \eps^{-1/q} \big)\, .$$

\subsection{Self-intersection local time: The upper bound}

Using the notation from the previous section, given $\eps>0$ we pick the integer~$n$ such that
$2^{-(n+1)q} \leq 2 \, C(q) \, \eps < 2^{-nq}$. Using the
super-additivity of $x \mapsto x^q$ for $q \geq 1$ we get
\[ \int_{-\infty}^\infty L^q(x,\sigma) \,  dx \geq \sum_{j=1}^{N(n)}
\int_{- \infty}^\infty L^q(x,\tau_{j-1}^{\ssup n}, \tau_{j}^{\ssup n}) \, dx \geq \eps \,
2^{-n+1} \sum_{j=1}^{2^n} Y_j \, . \] Hence, we get
\[ \p\Big\{ \int_{-\infty}^\infty L^q(x,\sigma) \, dx < \eps \Big\} \leq \p \Big\{
2^{1-n} \sum_{j=1}^{2^n} Y_j < 1 \Big\} = \p \{ S(2^n) > 0 \} \, ,
\] where $S(k) := \sum_{j=1}^k X_j$ for $X_j:=\frac 12 - Y_j$. By the simple large deviation bound 
for the sum of bounded random variables with negative mean, given 
in Section~\ref{boettcherup}, we deduce the existence of a constant $0<\varphi<1$ such that
\[ - \log \p\Big\{ \int_{-\infty}^\infty L^q(x,\sigma) \,  dx < \eps \Big\} 
\geq - \log \prob\{ S(2^n) \ge 0 \} 
\geq (-\log \varphi)2^n \geq \tilde{c}(q) \eps^{-\frac{1}{q}} \, , \] 
for the constant $\tilde{c}(q):= (-\log\varphi)(2^{-1 - 1/q}C(q)^{-1/q})> 0$.



\section{Outlook to future research}

Small value probabilities for intersection local times of Brownian motions in dimensions two and three
are considerably more difficult to handle, but in principle our method still applies. 
An analogue of Theorem~\ref{main2} for
Brownian motions in dimensions two and three is proved using the branching tree heuristic 
in~\cite{MS07}, see also~\cite{KM05} for partial results and their applications in multifractal analysis.
\smallskip

There is no direct analogue to Theorem~\ref{main3} for a higher dimensional Brownian motion. However,
our main results have natural analogues for random walks and in the random walk setting problems analogous 
to  Theorem~\ref{main3} can also be tackled in higher dimensions. This research project, together with some 
applications to weakly self-avoiding walks, is currently ongoing.
\smallskip

Finally, it is a natural question to ask whether the main results of the present paper can be extended 
from Brownian motion to L\'evy processes. It appears that the approach presented here may be
suited for such an extension, and further investigations in this problem are promising. \smallskip
\bigskip

\pagebreak[3]

{\bf Acknowledgments:} We thank the Nuffield Foundation for awarding an Undergraduate Research Bursary,
which allowed us to study the first example. The first author would like to thank Greg Lawler and Wenbo Li
for useful discussions and the EPSRC for support through grant~EP/C500229/1 
and an Advanced Research Fellowship.
\bigskip

\bibliographystyle{chicago}
\bibliography{small}

\end{document}